\documentclass[11pt,reqno]{amsart}
\usepackage{hyperref}
\usepackage{amssymb,amsmath,amsfonts,latexsym}
\usepackage{bm}
\usepackage{mathtools}
\usepackage{enumerate}
\usepackage{array,graphics,color}
\allowdisplaybreaks

\numberwithin{equation}{section}
\newtheorem{dfn}{Definition}[section]
\newtheorem{thm}[dfn]{Theorem}

\newtheorem{rmrk}[dfn]{Remark}

\DeclarePairedDelimiterX{\norm}[1]{\lVert}{\rVert}{#1}
\DeclarePairedDelimiterX{\bnorm}[1]{\big\lVert}{\big\rVert}{#1}
\DeclarePairedDelimiterX{\Bnorm}[1]{\Big\lVert}{\Big\rVert}{#1}

\newcommand{\exch}{\overset{\textbf{\textit{(Intc)}}}{=\joinrel=}}
\newcommand{\adj}{\overset{\textbf{\textit{(Adj)}}}{=\joinrel=}}
\newcommand{\sgn}{\overset{\textbf{\textit{(Sgn)}}}{=\joinrel=}}

\newcommand{\Z}{\mathbb{Z}}
\newcommand{\T}{\mathbb{T}}
\newcommand{\N}{\mathbb{N}}

\newcommand{\C}{\mathbb{C}}

\newcommand{\hdct}{{H}^2_{\mathbb{C}^2}(\mathbb{D})}

\newcommand{\hdcc}{{H}^2_{\mathbb{C}}(\mathbb{D})}

\begin{document}
	
	\title[Characterization of C-Symmetric Toeplitz operators]{ Characterization of C-Symmetric Toeplitz operators for a Class of Conjugations in Hardy Spaces}
	
	\author[Chattopadhyay] {Arup Chattopadhyay}
	\address{Department of Mathematics, Indian Institute of Technology Guwahati, Guwahati, 781039, India}
	\email{arupchatt@iitg.ac.in, 2003arupchattopadhyay@gmail.com}

	\author[Das]{Soma Das}
	\address{Department of Mathematics, Indian Institute of Technology Guwahati, Guwahati, 781039, India}
	\email{soma18@iitg.ac.in, dsoma994@gmail.com}
	
	\author[Pradhan]{Chandan Pradhan}
	\address{Department of Mathematics, Indian Institute of Technology Guwahati, Guwahati, 781039, India}
	\email{chandan.math@iitg.ac.in, chandan.pradhan2108@gmail.com}
	
	\author[Sarkar] {Srijan Sarkar}
	\address{Department of Mathematics, Indian Institute of Science, Bangalore 560012, India}
	\email{srijans@iisc.ac.in, srijansarkar@gmail.com}

	
	
	\begin{abstract}
	In this article, we introduce a new class of conjugations in the scalar valued Hardy space $H^2_{\mathbb{C}}(\mathbb{D})$ and provide a characterization of a complex symmetric Toeplitz operator $T_{\phi}$ with respect to these newly introduced conjugations in various cases. Moreover, we obtain a characterization of a complex symmetric block Toeplitz operator $T_{\Phi}$ on the vector valued Hardy space $\hdct$ with respect to certain conjugations introduced in \cite{CamaGP,KangKoLee,LeeKo2019}.
	\end{abstract}
	\maketitle
	
	\section{Introduction and Preliminaries}
	Complex symmetric operators on Hilbert spaces are natural generalizations of complex symmetric matrices, and the study of complex symmetric (in short C-symmetric) operators was initiated by Garcia, Putinar and Wogen in \cite{GP06,GP07,GW09,GW10}. A well known class of operators, namely all normal operators, Hankel operators and truncated Toeplitz operators are included in the class of complex symmetric operators. For more on complex symmetric operators and related stuff including historical comments we refer the reader to \cite{CamaGP,G06,GP06,GP07,GPP14,GW09,GW10,KangKoLee,LeeKo2016,LeeKo2019} and the references cited therein. 
	
	The following concept is a straightforward generalization of the conjugate-linear map  $z \to \overline{z}$ on the one-dimensional Hilbert space $\mathbb{C}$. 
 \begin{dfn}\cite{GP07}
		A conjugation on a complex Hilbert space $\mathcal{H}$ is a function $ C:\mathcal{H} \to \mathcal{H} $ which satisfies the following three properties:
		\begin{enumerate}[(i)]
			\item conjugate linear: $C(\alpha x+\beta y)=\overline{\alpha}Cx +\overline{\beta}Cy$, for all $x,y \in \mathcal{H}$ and $\forall \alpha , \beta \in \C$,
			\item  involutive: $C^2=I$, 
			\item isometric: $\norm {Cx}= \norm {x}.$ 
		\end{enumerate}
	\end{dfn}
In this connection it is worth mentioning that, Garcia and Putinar have shown in \cite{GP06} that for any given conjugation $C$ on a separable Hilbert space $\mathcal{H}$ there exists an orthonormal basis $\{e_n:n\in \N _0 \}$ such that $Ce_n=e_n$, where $\N _0$ denotes the set of all non-negative integers. Let $\mathcal{H}$ be a separable complex Hilbert space and  $\mathcal{B}(\mathcal{H})$ denote the set of all bounded linear operators on $\mathcal{H}$. 
	\begin{dfn}
		An operator $T\in \mathcal{B}(\mathcal{H})$ is called  $C$-symmetric if there exists a conjugation $C$ on $\mathcal{H}$  such that  $CTC=T^*$. If $T$ is $C$-symmetric for some conjugation $C$, then $T$ is called complex symmetric.
	\end{dfn}
Let $\mathbb{D}=\{z\in \mathbb{C}:~|z|<1\}$ denote the unit disc in the complex plane and $\mathbb{T}=\{z\in \mathbb{C}:~|z|=1\}$ be the unit circle. Let $L^2(\mathbb{T})$ be the Lebesgue (Hilbert) space on $\mathbb{T}$ and let $L^{\infty}(\mathbb{T})$ be the Banach
space of all essentially bounded functions on $\mathbb{T}$. Now it is well-known that $\{e_n(z) = z^n:~n\in \mathbb{Z}\}$ is an orthonormal basis for $L^2(\mathbb{T})$, where $\mathbb{Z}$ is the set of all integers. Therefore, if $f\in L^2(\mathbb{T})$, then the function $f$ can be expressed
as $$f(z)=\sum_{n=-\infty}^{\infty}\hat{f}(n)z^n,$$ where $\hat{f}(n)$ denotes the $n$-th Fourier coefficient of $f$ and $\sum\limits_{n=-\infty}^{\infty}|\hat{f}(n)|^2 <\infty$. 	Recall that $\mathcal{H}$-valued Hardy space over the unit disc $\mathbb{D}$ in $\mathbb{C}$ is denoted by $	H^2_{\mathcal{H}}(\mathbb{D})$ and defined by 
\begin{equation*}
	H^2_{\mathcal{H}}(\mathbb{D}):=\Big\{f(z)=\sum_{n=0}^{\infty}a_nz^n:~\|f\|_{H^2_{\mathcal{H}}(\mathbb{D})}^2:=\sum_{n=0}^{\infty}\|a_n\|_{\mathcal{H}}^2,~z\in \mathbb{D},~a_n\in \mathcal{H}\Big\}.
\end{equation*}
In this article, we mainly focus on two particular vector valued Hardy spaces, namely the classical Hardy spaces $H^2_{\mathbb{C}}(\mathbb{D})$ and $H^2_{\mathbb{C}^2}(\mathbb{D})$ corresponding to $\mathcal{H}=\mathbb{C}$ and $\mathcal{H}=\mathbb{C}^2$ respectively. For any $\phi\in L^{\infty}(\mathbb{T})$, the Toeplitz operator $T_{\phi}:H^2_{\mathbb{C}}(\mathbb{D})\longrightarrow H^2_{\mathbb{C}}(\mathbb{D})$ is defined by the formula 
\[
T_{\phi}(f)=P(\phi f),\quad f\in H^2_{\mathbb{C}}(\mathbb{D}),
\]
where $P$ denotes the orthogonal projection of $L^2(\mathbb{T})$ onto $H^2_{\mathbb{C}}(\mathbb{D})$. It is well known that $T_{\phi}$ is bounded if and only if $\phi\in L^{\infty}(\mathbb{T})$, and moreover, $\|T_{\phi}\|=\|\phi\|_{\infty}$. Note that we can also identify the Hardy space $H^2_{\mathcal{H}}(\mathbb{D})$ as a Hilbert space tensor product between $H^2_{\mathbb{C}}(\mathbb{D})$ and $\mathcal{H}$, that is $H^2_{\mathcal{H}}(\mathbb{D})=H^2_{\mathbb{C}}(\mathbb{D})\otimes \mathcal{H}$. 

Let $L^2_{\mathbb{C}^2}(\mathbb{T})=L^2(\mathbb{T})\otimes \mathbb{C}^2$, and let $L^{\infty}_{M_2}(\mathbb{T})=L^{\infty}(\mathbb{T})\otimes M_2$, where $M_2$ is the set of all $2\times 2$ complex matrices. Now for $\Phi\in L^{\infty}_{M_2}(\mathbb{T})$, the block Toeplitz operator with symbol $\Phi$
is the operator $T_{\Phi}$ on the vector-valued Hardy space $H^2_{\mathbb{C}^2}(\mathbb{D})$ defined by 
\[
T_{\Phi}(f)=\widetilde{P}(\Phi f),\quad f\in H^2_{\mathbb{C}^2}(\mathbb{D}),
\]
where $\widetilde{P}$ is the orthogonal projection of $L^2_{\mathbb{C}^2}(\mathbb{T})$ onto $H^2_{\mathbb{C}^2}(\mathbb{D})$. In particular, if
$	\Phi= \begin{bmatrix}
		\phi_1 & \phi_2 \\
		\phi_3 & \phi_4
	\end{bmatrix}$
where $\phi_i\in L^{\infty}(\mathbb{T})$ for $1\leq i\leq 4$, then the block Toeplitz operator has the
following representation:
\begin{equation*}
	T_{\Phi}= \begin{bmatrix}
		T_{\phi_1} & T_{\phi_2} \\
		T_{\phi_3} & T_{\phi_4}
	\end{bmatrix}.
\end{equation*}
For more on block Toeplitz operator and related topics we refer the reader to \cite{CurtHL}.

The study of complex symmetric operators, Toeplitz operators and block Toeplitz operators provides important connections with various problems in the ﬁeld of physics and most  importantly in the field of mechanics \cite{BE07,BFGJ,GPP14,Prunel}. Normal operators are examples of complex symmetric operators and the characterization of
normal Toeplitz operators was given by Brown and Halmos in \cite{BH}. In other words, they proved that $T_{\phi}$ is normal if and only if $\phi=\alpha + \beta \rho$ for some real-valued
function $\rho \in L^{\infty}(\mathbb{T})$ and $\alpha, \beta \in \mathbb{C}$. Note that, if $\phi\in L^{\infty}(\mathbb{T})$, then $T_{\phi}$ may not be a complex symmetric operator. Therefore, in general, it is a difficult problem to describe when a Toeplitz operator is complex symmetric. In this direction, recently, K. Guo and S. Zhu \cite{GZ} have raised the following interesting question: {\emph{Characterize a complex symmetric Toeplitz operator on the Hardy space $H^2_{\mathbb{C}}(\mathbb{D})$}.  This question has motivated researchers to identify special classes of conjugations on Hardy spaces}. More precisely, for certain conjugations $C$ with explicit forms, it is an interesting question to characterize $C$-symmetric Toeplitz operators. Recently, E. Ko and J.E. Lee in \cite{LeeKo2016} gave a characterization of a complex symmetric Toeplitz operator $T_{\phi}$ on $H^2_{\mathbb{C}}(\mathbb{D})$ with respect to some special conjugations. More precisely, they considered the family of conjugations $C_{\mu,\lambda}$ on $H^2_{\mathbb{C}}(\mathbb{D})$ deﬁned by
\[
C_{\mu,\lambda}f(z) = \mu \overline{f(\lambda \bar{z})}
\]
for $\mu,\lambda\in \mathbb{T}$, and proved the following theorem:
\begin{thm}
 If $\phi\in L^{\infty}(\mathbb{T})$, then $T_{\phi}$ is $C_{\mu,\lambda}$-symmetric if and only if $\hat{\phi}(-n)=\hat{\phi}(n)\lambda^n$ for all $n\in \mathbb{Z}$, where $\hat{\phi}(n)$ is the $n$-th Fourier coefficients of $\phi$. 
\end{thm}
In this context, a recent result by S. Waleed Noor \cite{Noor} proves the following result: if $T_{\phi}$ is complex symmetric on $H^2_{\mathbb{C}}(\mathbb{D})$ with continuous symbol $\phi$ on $\mathbb{T}$, 
then $\phi(\mathbb{T})$ is a nowhere winding curve. We refer the reader to articles \cite{LLC,LYL} for important results in the study of complex symmetric Toeplitz operators on Bergman spaces and Dirichlet spaces. Moreover, very recently D. Kang, E. Ko, and J.E. Lee provide a characterization of complex symmetric block Toeplitz operator $T_{\Phi}$ with respect to some special conjugations on the vector-valued Hardy space $H^2_{\mathbb{C}^2}(\mathbb{D})$.

Motivated by all these works (most importantly, \cite{ KangKoLee,LeeKo2016, LeeKo2019}), our principle aim in this article is to give some characterizations of complex symmetric  Toeplitz operators $T_{\phi}$  and $T_{\Phi}$ on the Hardy spaces $H^2_{\mathbb{C}}(\mathbb{D})$ and $H^2_{\mathbb{C}^2}(\mathbb{D})$, respectively with respect to certain new conjugations defined as follows: Let $p\in \mathbb{N}$ and let $S_p$ denote the symmetric group defined over a finite set of $p$ symbols consisting of the permutations that can be performed on the $p$ symbols. For  $\sigma\in S_p$, we denote $O(\sigma)$ as the order of the permutation $\sigma$. Now for $\sigma\in S_p$ with $O(\sigma)=2$, let $C_\sigma:\hdcc \mapsto \hdcc$ be defined by
\begin{align}\label{mainconju}
	C_\sigma \left(\sum_{k=0}^{\infty}\sum_{m=0}^{p-1}a_{m+pk}z^{m+pk}\right)= \sum_{k=0}^{\infty}\sum_{m=0}^{p-1}\overline{\sigma(a_{m+pk})}z^{m+pk},
\end{align}
where for fixed $p$ and $k$, $\sigma$ is a permutation on the set $\{a_{pk},a_{1+pk},\ldots,a_{(p-1)+pk}\}$. Then it is easy to verify from the definition that $C_{\sigma}$ is a conjugation on $\hdcc$. 

This paper is organized as follows: In section 2, we provide a characterization of complex symmetric Toeplitz operators $T_{\phi}$ with respect to a special case of \eqref{mainconju}, that is, the conjugation  $C_p^{i,j}$ for some fix $p\in \mathbb{N}$ and $i,j\in \mathbb{N}$ such that $i\neq j$ (see Theorem~\ref{mainthm}). Section 3 deals with the characterization of Toeplitz operators $T_{\phi}$ with respect to a conjugation $C_n$ on $H^2_{\mathbb{C}}(\mathbb{D})$ which is again a special case of \eqref{mainconju}. In section 4, we give a characterization of block Toeplitz operators $T_{\Phi}$ with respect to the conjugations $C$ (see \eqref{conjuvec1} ) and $\widetilde{C}$ (see \eqref{conjuvec2}) on $H^2_{\mathbb{C}^2}(\mathbb{D})$, respectively (see Theorem~\ref{thm3}, Theorem~\ref{thm4} and Theorem~\ref{thm5}) that were introduced earlier in \cite{CamaGP,KangKoLee,LeeKo2019}.

\section{Transpositions type of  Conjugations}
In this section we consider a class of conjugations which are special case of \eqref{mainconju} and study the complex symmetry of the Topelitz operator $T_{\phi}$ on $H^2_{\mathbb{C}}(\mathbb{D})$ with respect to those conjugations. 
In other words, for fix $p\in \N$, we choose $i,j\in \mathbb{N}$ such that $0\leq i < j < p$ and define the map $C_{p}^{i,j}: \hdcc \mapsto \hdcc$ by
\begin{align}\label{splconju}
C_{p}^{i,j}\left(\sum_{k=0}^{\infty}a_kz^k\right) \mapsto \sum_{k=0}^{\infty}\bar{a}_{j+pk}z^{i+pk}  +\sum_{k=0}^{\infty}\bar{a}_{i+pk}z^{j+pk}+\sum_{k=0}^{\infty}\sum_{\underset{m\neq i,j}{m=0}}^{p-1}\bar{a}_{m+pk}z^{m+pk}.
\end{align}
 Then it follows from the definition above that $C_{p}^{i,j}$ is a conjugation on $\hdcc$. Our next aim is to investigate some necessary and sufficient conditions on the symbol $\phi\in \L^{\infty}(\mathbb{T})$ for which the Toeplitz operator $T_{\phi}$ is complex symmetric with respect to the conjugation $C_p^{i,j}$.  To that end, let  $\phi (z)=\sum\limits_{k=-\infty}^{\infty}\hat{\phi}(k)z^k \in L^\infty (\T)$ and consider the associated Toeplitz operator $T_\phi$. Next we assume that $T_\phi$ is complex symmetric with respect to the conjugation $C_{p}^{i,j}$, that is
\begin{align}\label{3}
T_\phi C_{p}^{i,j} =C_{p}^{i,j} T_{\bar{\phi}}.
\end{align}
Now it is well known that the set $\{z^{r+pk}:k\geq 0, 0\leq r\leq p-1 \}$ forms an orthonormal basis of $\hdcc$ and therefore corresponding to $r=j$ we have 
\begin{align}\label{eq0}
T_\phi C_{p}^{i,j}(z^{j+pk}) = T_\phi (z^{i+pk})
=\sum_{m=0}^{\infty}\hat{\phi}\big(m-i-pk\big)z^m,
\end{align}
and
\begin{align}\label{eq1}
C_{p}^{i,j}T_{\bar{\phi}} (z^{j+pk}) 
\nonumber & = C_{p}^{i,j}\left(\sum_{m=0}^{\infty}\overline{\hat{\phi}\big(-(m-j-pk)\big)}z^m \right)\\
\nonumber &= \sum_{m=0}^{\infty}\hat{\phi}\big(-(i+pm-j-pk)\big)z^{j+pm} + \sum_{m=0}^{\infty}\hat{\phi}\big(-(j+pm-j-pk)\big)z^{i+pm}\\
&\hspace{1in} +\sum_{m=0}^{\infty}\sum_{\underset{t\neq i,j}{t=0}}^{p-1}\hat{\phi}\big(-(t+pm-j-pk)\big)z^{t+pm}.
\end{align}
Thus by substituting $m-k=l \in \Z $  in \eqref{eq0} and \eqref{eq1}, and using \eqref{3} we get
\begin{align}\label{eq2}
\hat{\phi}(pl)= \hat{\phi}(-pl),~ \hat{\phi}(j-i+pl)= \hat{\phi}(j-i-pl),~ \hat{\phi}(t-i+pl)= \hat{\phi}(j-t-pl).
\end{align}
Moreover, for $r=i$, by repeating the above calculations and using the relation $T_\phi C_{p}^{i,j}(z^{i+pk}) =C_{p}^{i,j} T_{\bar{\phi}}(z^{i+pk})$ we conclude for $l\in \Z$ that
\begin{align}\label{eq3}
\hat{\phi}(i-j+pl)= \hat{\phi}(i-j-pl), ~\hat{\phi}(pl)= \hat{\phi}(-pl), ~ \hat{\phi}(t-j+pl)= \hat{\phi}(i-t-pl).
\end{align}
Now for $0\leq r\leq p-1$ such that $r\neq i,j$ we have
\begin{align*}
T_\phi C_{p}^{i,j}(z^{r+pk})= \sum_{m=0}^{\infty}\hat{\phi}\big(m-r-pk\big)z^m,
\end{align*}
and 
\begin{align*}
C_{p}^{i,j}T_{\bar{\phi}} (z^{r+pk}) & = C_{p}^{i,j}\Big(\sum_{m=0}^{\infty}\overline{\hat{\phi}\big(-(m-r-pk)\big)}z^m \Big)\\
& =\sum_{m=0}^{\infty}\hat{\phi}\big(-(j+pm-r-pk)\big)z^{i+pm}+ \sum_{m=0}^{\infty}\hat{\phi}\big(-(i+pm-r-pk)\big)z^{j+pm} \\
&\hspace{0.5in}+\sum_{m=0}^{\infty}\hat{\phi}\big(-(r+pm-r-pk)\big)z^{r+pm}+\sum_{m=0}^{\infty}\sum_{\underset{s\neq r,i,j}{s=0}}^{p-1}\hat{\phi}\big(-(s+pm-r-pk)\big)z^{s+pm}.
\end{align*}
Finally, by using \eqref{3} we obtain the following conditions
 \begin{align}\label{eq4}
\nonumber \hat{\phi}(i-r+pl)= \hat{\phi}(r-j-pl), ~&~ \hat{\phi}(j-r+pl)= \hat{\phi}(r-i-pl),\\
 \hat{\phi}(s-r+pl)= \hat{\phi}(r-s-pl), ~& ~\hat{\phi}(pl)= \hat{\phi}(-pl),
 \end{align}
where $m-k=l \in \Z$ and $0\leq s\leq p-1$ such that $s\neq i,j,r$.
All the conditions obtained in \eqref{eq2},\eqref{eq3} and \eqref{eq4} are enlisted in the following table
\begin{align}\label{eq5}
	\begin{tabular}{ |c| } 
		\hline
		If $l\in \Z,0\leq a\leq p-1$ such that $a\neq i,j$, and $0\leq b\leq p-1$ such that $b\neq a,i,j$, then   \\ 
		\hline
		 $\hat{\phi}\Big(pl\Big)= \hat{\phi}\Big(-pl\Big)$,\\
		 \hline 
		 $\hat{\phi}\Big(j-i+pl\Big)= \hat{\phi}\Big(j-i-pl\Big)$ , $\hat{\phi}\Big(i-j+pl\Big)= \hat{\phi}\Big(i-j-pl\Big)$,\\
		 \hline
		$\hat{\phi}\Big(a-i+pl\Big)= \hat{\phi}\Big(j-a-pl\Big)$ , $\hat{\phi}\Big(a-j+pl\Big)= \hat{\phi}\Big(i-a-pl\Big)$,\\
		\hline
		$\hat{\phi}\Big(i-a+pl\Big)= \hat{\phi}\Big(a-j-pl\Big)$ , $\hat{\phi}\Big(j-a+pl\Big)= \hat{\phi}\Big(a-i-pl\Big)$,\\
		\hline
		$\hat{\phi}\Big(b-a+pl\Big)= \hat{\phi}\Big(a-b-pl\Big)$.\\
		\hline
	\end{tabular}
\end{align}
Our next aim is to simplify the above relations between the Fourier coefficients of $\phi$ obtained in \eqref{eq5} by assuming some more restrictions on $i,j$ and $p$.
\vspace{0.1in}

\textbf{Case I:} First we assume that $p$ is even and $\vert i-j\vert =\frac{p}{2}$. That is, $j-i=\frac{p}{2}$ and $i-j=-\frac{p}{2}.$ Therefore, $j=\frac{p}{2} +i$
and hence the above conditions mentioned in \eqref{eq5} becomes
\begin{align}\label{eq6}
	\begin{tabular}{ |c| } 
		\hline
		If $l\in \Z,0\leq a\leq p-1$ such that $a\neq i,j$, and $0\leq b\leq p-1$ such that $b\neq a,i,j$, then   \\ 
		\hline
		$\hat{\phi}\Big(pl\Big)= \hat{\phi}\Big(-pl\Big)$,\\
		\hline 
		$\hat{\phi}\Big(\frac{p}{2}+pl\Big)= \hat{\phi}\Big(\frac{p}{2}-pl\Big)$ , $\hat{\phi}\Big(-\frac{p}{2}+pl\Big)= \hat{\phi}\Big(-\frac{p}{2}-pl\Big)$,\\
		\hline
		$\hat{\phi}\Big(a-i+pl\Big)= \hat{\phi}\Big(\frac{p}{2}+i-a-pl\Big)$ , $\hat{\phi}\Big(a-\frac{p}{2}-i+pl\Big)= \hat{\phi}\Big(i-a-pl\Big)$,\\
		\hline
		$\hat{\phi}\Big(i-a+pl\Big)= \hat{\phi}\Big(a-\frac{p}{2}-i-pl\Big)$ , $\hat{\phi}\Big(\frac{p}{2}+i-a+pl\Big)= \hat{\phi}\Big(a-i-pl\Big)$,\\
		\hline
		$\hat{\phi}\Big(b-a+pl\Big)= \hat{\phi}\Big(a-b-pl\Big)$.\\
		\hline
	\end{tabular}
\end{align} 
The following two tables consist of different values of $a$ that are essential in the sequel.
\begin{align}\label{eq7}
	\begin{tabular}{ |c| c| c| c| c| c| c| c| c|} 
		\hline
		$a$ & $0$ & $1$ & $\cdots$ & $i-1$ & $i+1$& $i+2$ & $\cdots$ & $\frac{p}{2}+i-1$\\[2ex]
		\hline
		${a-i}$ & $-i$ & $1-i$  & $\cdots$ & $-1$ & $1$ & $2$ & $\cdots$ & $\frac{p}{2}-1$ \\[2ex]
		\hline
		${i-a}$ & $i$ & $i-1$ &  $\cdots$ & $1$ & $-1$ & $-2$ & $\cdots$ & $1-\frac{p}{2}$ \\[2ex]
		\hline
		${\frac{p}{2}+i-a}$ & $\frac{p}{2}+i$ & $\frac{p}{2}+i-1$ &  $\cdots$ & $\frac{p}{2}+1$ & $\frac{p}{2}-1$ & $\frac{p}{2}-2$ & $\cdots$ & $1$ \\[2ex]
		\hline
		${a-\frac{p}{2}-i}$ & $-\frac{p}{2}-i$ & $1-\frac{p}{2}-i$ & $\cdots$ & $-\frac{p}{2}-1$ & $1-\frac{p}{2}$ & $2-\frac{p}{2}$ & $\cdots$ & $-1$ \\[2ex]
		\hline	
	\end{tabular}
\end{align}

\begin{align}\label{eq8}
		\begin{tabular}{|c |c| c| c| c| c|} 
			\hline
			$a$ & $\frac{p}{2}+i+1$ & $\frac{p}{2}+i+2$ & $\cdots$ & $p-2$ & $p-1$ \\[2ex]
			\hline
			$a-i$ &$\frac{p}{2}+1$ & $\frac{p}{2}+2$ & $\cdots$ & $p-2-i$ & $p-1-i$ \\[2ex]
			\hline
			$i-a$ &$-\frac{p}{2}-1$ & $-\frac{p}{2}-2$ & $\cdots$ & $i-p+2$ & $i-p+1$ \\[2ex]
			\hline
			$\frac{p}{2}+i-a$ &$-1$ & $-2$ & $\cdots$ & $-\frac{p}{2}+i+2$ & $-\frac{p}{2}+i+1$ \\[2ex]
			\hline
			$a-\frac{p}{2}-i$ &$1$ & $2$ & $\cdots$ & $\frac{p}{2}-i-2$ & $\frac{p}{2}-i-1$ \\[2ex]
			\hline	
		\end{tabular}
\end{align}

\textbf{Sub-case I:} Suppose {\bm{$\frac{p}{2}$}} is even.
Then we have the following $\frac{p}{4}$ many pairs:
\begin{equation}\label{eq12}
\Big(1,\frac{p}{2}-1 \Big),\Big(2,\frac{p}{2}-2 \Big),\cdots,\Big(\frac{p}{4}-1,\frac{p}{4}+1 \Big),\Big(\frac{p}{4},\frac{p}{4} \Big).	
\end{equation}
Proceed further, we need the following conditions as mentioned in the 4th row of the table \eqref{eq6}, that is for $l\in \mathbb{Z}$, and for $0\leq a\leq p-1$ such that $a\neq i,\frac{p}{2}+i$, we have
\begin{align}\label{eq9}
\hat{\phi}\Big(a-i+pl\Big)= \hat{\phi}\Big(\frac{p}{2}+i-a-pl\Big),~ \hat{\phi}\Big(a-\frac{p}{2}-i+pl\Big)= \hat{\phi}\Big(i-a-pl\Big). 
\end{align}
Therefore using the column corresponding to $a=i+1$ in \eqref{eq7} and using \eqref{eq9} we get
\begin{align*}
\hat{\phi}\left(1+pl \right)  = \hat{\phi}\left(\frac{p}{2}-1-pl \right)
= \hat{\phi}\left(-(\frac{p}{2}+1)-p(l-1) \right)
= \hat{\phi}\left(1+p(l-1) \right) \quad \forall \quad l\in \mathbb{Z},
\end{align*}
where at the last equality we have used the column corresponding to $a=\frac{p}{2}+i+1$ in \eqref{eq8}. Again, using \eqref{eq9} and using the columns corresponding to $a=i+2$ and $a=\frac{p}{2}+i+2$ in \eqref{eq7} and \eqref{eq8} respectively we obtain
\begin{align*}
\hat{\phi}\left(2+pl \right) = \hat{\phi}\left(\frac{p}{2}-2-pl \right)
= \hat{\phi}\left(-(\frac{p}{2}+2)-p(l-1) \right)
= \hat{\phi}\left(2+p(l-1) \right) \quad \forall \quad l\in \mathbb{Z}.
\end{align*}
Therefore by repeating the same argument as above and using \eqref{eq7}, \eqref{eq8} and \eqref{eq9} we conclude
\begin{align}\label{eq10}
	\hat{\phi}\left(r+pl \right) = \hat{\phi}\left(r+p(l-1) \right) \quad \forall \quad l\in \mathbb{Z}\quad \text{and}\quad 1\leq r\leq \frac{p}{4}.
\end{align}
Since $\phi \in L^\infty(\T) \subset L^2(\T)$ and hence $\sum\limits_{k=-\infty}^{\infty}\vert{\hat{\phi}(k)}\vert ^2 < \infty$, then the equation \eqref{eq10} yields 
\begin{align}\label{eq11}
	\hat{\phi}\left(r+pl \right)=0 \quad \forall \quad l\in \mathbb{Z} \quad  \text{and} \quad 1\leq r\leq \frac{p}{4}.
	\end{align}
Now by observing the symmetricity of the pair in \eqref{eq12} and using equations \eqref{eq9} and \eqref{eq11}  we conclude
\begin{align}\label{eq17}
\hat{\phi}\left(\frac{p}{4}+s+pl \right)= \hat{\phi}\left(\frac{p}{4}-s-pl \right)=0 \quad \forall \quad l\in \mathbb{Z}\quad \text{and}\quad 1\leq s\leq \frac{p}{4}-1,
\end{align}
Furthermore, using the 3rd row of the table \eqref{eq6} we get
\begin{align*}
\hat{\phi}\left(\frac{p}{2}+pl \right) & = \hat{\phi}\left(\frac{p}{2}-pl \right)
= \hat{\phi}\left(\frac{p}{2}-p-p(l-1) \right)\\
& = \hat{\phi}\left(-\frac{p}{2}-p(l-1) \right)
= \hat{\phi}\left(-\frac{p}{2}+p(l-1) \right)
= \hat{\phi}\left(\frac{p}{2}+p(l-2) \right), \quad \forall \quad l\in \Z,
\end{align*}
and hence by the similar argument as in \eqref{eq11} we conclude 
\begin{equation}\label{eq12}
	\hat{\phi}\left(\frac{p}{2}+pl \right)=0\quad \forall \quad l\in \Z.
\end{equation}
Next, we need the following conditions as mentioned in the 5th row of the table \eqref{eq6}, that is for $l\in \mathbb{Z}$, and for $0\leq a\leq p-1$ such that $a\neq i,\frac{p}{2}+i$, we have
\begin{align}\label{eq13}
\hat{\phi}\Big(i-a+pl\Big)= \hat{\phi}\Big(a-\frac{p}{2}-i-pl\Big), \hat{\phi}\Big(\frac{p}{2}+i-a+pl\Big)= \hat{\phi}\Big(a-i-pl\Big).
\end{align}
Again, by using the columns corresponding to $a=i+1, i+2$ in \eqref{eq7} and using \eqref{eq13} we get
\begin{align*}
\hat{\phi}\left(-1+pl \right) = \hat{\phi}\left(-\frac{p}{2}+1-pl \right)
= \hat{\phi}\left(\frac{p}{2}+1-p(l+1) \right)
= \hat{\phi}\left(-1+p(l+1) \right) \quad \forall \quad l\in \mathbb{Z},
\end{align*}
\begin{align*}
\hat{\phi}\left(-2+pl \right) = \hat{\phi}\left(-\frac{p}{2}+2-pl \right)
= \hat{\phi}\left(\frac{p}{2}+2-p(l+1) \right)
= \hat{\phi}\left(-2+p(l+1) \right) \quad \forall \quad l\in \mathbb{Z}.
\end{align*}
Therefore by repeating the same argument as above and using \eqref{eq7}, \eqref{eq8} and \eqref{eq13} we conclude
\begin{align*}
	\hat{\phi}\left(-r+pl \right) = \hat{\phi}\left(-r+p(l+1) \right) \quad \forall \quad l\in \mathbb{Z}\quad \text{and}\quad 1\leq r\leq \frac{p}{4}.
\end{align*}
Thus by the similar argument as in \eqref{eq11} we conclude 
\begin{align}\label{eq14}
	\hat{\phi}\left(-r+pl \right)=0 \quad \forall \quad l\in \mathbb{Z} \quad  \text{and} \quad 1\leq r\leq \frac{p}{4}.
\end{align}
Consequently, by using equations \eqref{eq13} and \eqref{eq14} we get
\begin{align}\label{eq15}
\hat{\phi}\left(\frac{p}{2}+r+pl \right)= \hat{\phi}\left(-r-pl \right)=0 \quad \forall \quad l\in \mathbb{Z} \quad  \text{and} \quad 1\leq r\leq \frac{p}{4}.
\end{align}	
Furthermore, using \eqref{eq14} we also conclude
\begin{align}\label{eq16}
\hat{\phi}\left(\frac{p}{2}+\frac{p}{4}+s+pl \right) = \hat{\phi}\left(-(\frac{p}{4}-s) +p(l+1) \right) =0 \quad \forall \quad l\in \mathbb{Z} \quad  \text{and} \quad 1\leq s\leq \frac{p}{4}-1.
\end{align}
Thus, by combining all the conditions obtained in \eqref{eq11}, \eqref{eq17}, \eqref{eq12}, \eqref{eq15} and \eqref{eq16} we get 
\begin{align*}
\hat{\phi}(r+pl)=0\quad \forall \quad l\in \mathbb{Z} \quad  \text{and} \quad 1\leq r\leq p-1.
\end{align*}

\textbf{Sub-case II:} Suppose {\boldmath{$\frac{p}{2}$}} is odd. Then we have the following $\dfrac{p/2-1}{2}$ many pairs:
$$\Big(1,\frac{p}{2}-1 \Big),\Big(2,\frac{p}{2}-2 \Big),\cdots,\left( \dfrac{p/2-1}{2},\dfrac{p/2-1}{2}+1\right).$$
Therefore by proceeding with the similar arguments as in \textbf{Sub-case I} we conclude,
\begin{align*}
	\hat{\phi}(r+pl)=0\quad \forall \quad l\in \mathbb{Z} \quad  \text{and} \quad 1\leq r\leq p-1.
\end{align*}
\vspace{0.1in}

{\bf Case II:} Here we assume that $p=mq+1$ for some natural number $m\geq 2$, $i=q-1,$ and $j=p-1$.
Now by rewriting the relations obtained in the table \eqref{eq5}, we have the following 
\vspace{0.1in}

$\bullet~\textbf{Interchange Rule: ~(Intc)}$
\begin{align}
	\label{eq18}~&\hat{\phi}(c+pl)=\hat{\phi}(d-pl)  \text{ for }~
	\begin{cases}
		& |c|, |d|\in\{1,\ldots,p-1\}\setminus\{p-q\} \text{ such that } |c+d|=p-q,\\	
		&   c=d=\pm(p-q).
	\end{cases}
\end{align}
\vspace{0.1in}

$\bullet~\textbf{Sign Rule: ~(Sgn)}$
\begin{align}
 \label{eq19}&\hat{\phi}(c+pl)=\hat{\phi}(-c-pl) \text{ for } 
	\begin{cases}
		& c\in \{1,\ldots,p-3\} \text{ if } q=1,\\	
		& c\in \{1,\ldots,p-2\} \text{ if } q\geq 2.
	\end{cases}
\end{align} 
 Our next aim is to show that $\hat{\phi}(k+pl)=0$ for all $k\in\{1,2,\ldots,p-1\}$ and for all $l\in\mathbb{Z}$ by using the above two rules \eqref{eq18} and \eqref{eq19}.\\

Suppose $q=1$ and $p=m+1$ is odd (that is, $m$ is even). Now for $m=2$, we have
\begin{align}\label{eq20}
	\nonumber&\hat{\phi}(1+pl)=\hat{\phi}(1+3l)\exch\hat{\phi}(1-3l)\adj\hat{\phi}(-2-3(l-1))\\
	&\exch\hat{\phi}(-2+3(l-1))\adj\hat{\phi}(1+3(l-2)) \quad \forall \quad l\in \mathbb{Z},
\end{align}
where we have used equation \eqref{eq18} and the symbol ${\textbf{\textit{(Adj)}}}$ stands for the adjustment of the Fourier coefficients. Furthermore, by using equation \eqref{eq18} we have for $m>2$ ($m$ is even) that 
\begin{equation}\label{eq22}
	\begin{split}
		&\hat{\phi}(1+pl)\exch\hat{\phi}((p-2)-pl)\adj\hat{\phi}(-2-p(l-1))\exch\hat{\phi}((3-p)+p(l-1))\adj\hat{\phi}(3+p(l-2))\\
		&=\cdots\adj \hat{\phi}((p-2)+p(l-(p-3)))\exch\hat{\phi}(1-p(l-(p-3)))\adj\hat{\phi}((1-p)-p(l-(p-2)))\\
		&\exch\hat{\phi}((1-p)+p(l-(p-2)))\adj\hat{\phi}(1+p(l-(p-1))) \quad \forall \quad l\in \mathbb{Z}.
	\end{split}
\end{equation}
Now if $q=1$ and $p=m+1$ is even (that is, $m$ is odd), then again by using equations \eqref{eq18} and \eqref{eq19} we get
\begin{equation}\label{eq23}
	\begin{split}
		\hat{\phi}(1+pl)\exch&\hat{\phi}((p-2)-pl)\adj\hat{\phi}(-2-p(l-1))=\cdots\adj \hat{\phi}(-(p-2)-p(l-(p-3)))\\
		\exch&\hat{\phi}(-1+p(l-(p-3)))\adj\hat{\phi}((p-1)+p(l-(p-2)))\\
		\exch&\hat{\phi}((p-1)-p(l-(p-2)))\adj\hat{\phi}(-1-p(l-(p-1)))\\
		\sgn&\hat{\phi}(1+p(l-(p-1))) \quad \forall \quad l\in \mathbb{Z}.
	\end{split}
\end{equation}
Therefore by the similar argument 
as in \eqref{eq11}, equations \eqref{eq20},\eqref{eq22} and \eqref{eq23} yields that
\begin{align}\label{eq24}
	\hat{\phi}(k+pl)=0  \quad \forall \quad l\in \mathbb{Z}
	\quad \text{and} \quad k\in\{1,2,\ldots,p-1\}.
\end{align}
Next we consider $q\geq 2$, then by using \eqref{eq18} we get
\begin{align}\label{eq25}
	\nonumber\hat{\phi}(1+pl)&\exch\hat{\phi}((p-q)-1-pl)\adj\hat{\phi}(-(q+1)-p(l-1))\\
	\nonumber&\exch\hat{\phi}(-p+(2q+1)+p(l-1))\adj\hat{\phi}((2q+1)+p(l-2))\\
	\nonumber&\hspace{2in}\vdots\\
	&\adj\hat{\phi}\Big((-1)^{(m-1)}\big((m-1)q+1\big)+(-1)^{(m-1)}p\big(l-(m-1)\big)\Big) \quad \forall \quad l\in \mathbb{Z}.
\end{align}	
To proceed further, let us denote $l^{(k)}=l-k(m-1),\quad  k\in\mathbb{N}$,\quad  $l\in\mathbb{Z}$,\quad  $m\geq 2$.
\vspace{0.1in}

{\bf Sub-case I:} Suppose $p=mq+1$ such that $m$ is odd. Then the above equation \eqref{eq25} yields
\begin{align}\label{eq26}
	&\nonumber\hat{\phi}(1+pl)=\hat{\phi}\Big(\big((m-1)q+1\big)+p\big(l-(m-1)\big)\Big)=\hat{\phi}((p-q)+pl^{(1)})\exch\hat{\phi}((p-q)-pl^{(1)})\\
	\nonumber&\adj\hat{\phi}(-q-p(l^{(1)}-1))\exch\hat{\phi}((2q-p)+p(l^{(1)}-1))\adj\hat{\phi}(2q+p(l^{(1)}-2))\\
	\nonumber&\exch\cdots\adj\hat{\phi}(-mq-p(l^{(1)}-m))=\hat{\phi}\Big(-(p-1)-p\big(l^{(1)}-m\big) \Big)\\
	&\adj\hat{\phi}\Big(1-p\big(l^{(1)}-(m-1)\big) \Big)=\hat{\phi}\big(1-pl^{(2)} \big) \quad \forall \quad l\in \mathbb{Z},
\end{align}
where we have used equation \eqref{eq18}. Moreover, again by using the equation \eqref{eq18} we get
\begin{align}\label{eq27}
	\nonumber&\hat{\phi}\big(1-pl^{(2)} \big)\exch\hat{\phi}\Big((q-p)-1+pl^{(2)} \Big)\adj\hat{\phi}\Big((q-1)+p(l^{(2)}-1) \Big)\\
	\nonumber&\exch\hat{\phi}\Big(p-(2q-1)-p(l^{(2)}-1)\Big)\adj\hat{\phi}\Big(-(2q-1)-p(l^{(2)}-2)\Big)\\
	\nonumber&\exch\cdots=\adj\hat{\phi}\Big((-1)^{m-1}(mq-1)+(-1)^{m-1}p(l^{(2)}-m) \Big)=\hat{\phi}\Big(p-2+p\big(l^{(2)}-m\big) \Big)\\
	&\adj\hat{\phi}\Big(-2+p\big(l^{(2)}-(m-1)\big) \Big)=\hat{\phi}(-2+pl^{(3)}) \quad \forall \quad l\in \mathbb{Z}.
\end{align}	
Therefore by repeating the similar arguments as in \eqref{eq26} and \eqref{eq27} we conclude
\begin{align}\label{eq28}
	\nonumber\hat{\phi}(1+pl)=\hat{\phi}(1-pl^{(2)})=\hat{\phi}(-2+pl^{(3)})=\cdots&=\begin{cases}
		\hat{\phi}\Big(-(q-1)+pl^{(q)}\Big) &\text{ if $q$ is odd}\\
		\hat{\phi}\Big((q-1)-pl^{(q)}\Big) &\text{ if $q$ is even}
	\end{cases}\\
	\nonumber&\adj\begin{cases}
		\hat{\phi}\Big((p-q)+1+p\big(l^{(q)}-1\big)\Big) &\text{ if $q$ is odd}\\
		\hat{\phi}\Big((q-p)-1-p\big(l^{(q)}-1\big)\Big) &\text{ if $q$ is even}
	\end{cases}\\
	\nonumber&\exch\begin{cases}
		\hat{\phi}\Big(-1-p(l^{(q)}-1)\Big) &\text{ if $q$ is odd}\\
		\hat{\phi}\Big(1+p(l^{(q)}-1)\Big) &\text{ if $q$ is even}
	\end{cases}\\
	&\sgn\hat{\phi}\Big(1+p(l^{(q)}-1)\Big) \quad \forall \quad l\in \mathbb{Z}.
\end{align}
Therefore by combining all the conditions obtained in \eqref{eq26}, \eqref{eq27} and \eqref{eq28} we have the following chain of relations
\begin{equation}\label{eq29}
	\begin{split}
		&\hat{\phi}(1+pl)=\hat{\phi}(-(q+1)-p(l-1))=\hat{\phi}((2q+1)+p(l-2))=\cdots=\hat{\phi}\Big((p-q)+pl^{(1)}\Big)\\
		=&\hat{\phi}(-q+p(l^{(1)}-1))=\hat{\phi}(2q+p(l^{(1)}-2))=\cdots=\hat{\phi}(-mq-p(l^{(1)}-m))\\
		=&\hat{\phi}\big(1-pl^{(2)}\big)=\hat{\phi}(-2+pl^{(3)})=\cdots=
		\begin{cases}
			\hat{\phi}\Big(-(q-1)+pl^{(q)}\Big) &\text{ if $q$ is odd}\\
			\hat{\phi}\Big((q-1)-pl^{(q)}\Big) &\text{ if $q$ is even}
		\end{cases}\\
		=&\hat{\phi}\Big(1+p(l^{(q)}-1)\Big) \quad \forall \quad l\in \mathbb{Z}.
	\end{split}
\end{equation}
{\bf Sub-case II:}  Suppose $p=mq+1$ such that $m$ is even. Then by applying \eqref{eq18}, the equation \eqref{eq25} becomes
\begin{align}\label{eq30}
	\nonumber&\hat{\phi}(1+pl)=\hat{\phi}\Big(-\big((m-1)q+1\big)-p\big(l-(m-1)\big)\Big)=\hat{\phi}(-(p-q)-pl^{(1)})\exch\hat{\phi}(-(p-q)+pl^{(1)})\\
	\nonumber&\adj\hat{\phi}(q+p(l^{(1)}-1))\exch\hat{\phi}((p-2q)-p(l^{(1)}-1))\adj\hat{\phi}(-2q-p(l^{(1)}-2))\\\nonumber&\exch\cdots\adj\hat{\phi}(-mq-p(l^{(1)}-m))=\hat{\phi}\Big(-(p-1)-p\big(l^{(1)}-m\big) \Big)\adj\hat{\phi}\Big(1-p\big(l^{(1)}-(m-1)\big) \Big)\\
	\nonumber&=\hat{\phi}\big(1-pl^{(2)} \big)\exch\hat{\phi}\Big((q-p)-1+pl^{(2)} \Big)\adj\hat{\phi}\Big((q-1)+p(l^{(2)}-1) \Big)\\
	\nonumber&\exch\hat{\phi}\Big(p-(2q-1)-p(l^{(2)}-1)\Big)\adj\hat{\phi}\Big(-(2q-1)-p(l^{(2)}-2)\Big)\\
	\nonumber&\exch\cdots\adj\hat{\phi}\Big((-1)^{m-1}(mq-1)+(-1)^{m-1}p(l^{(2)}-m) \Big)=\hat{\phi}\Big(2-p-p\big(l^{(2)}-m\big) \Big)\\
	\nonumber&\adj\hat{\phi}\Big(2-p\big(l^{(2)}-(m-1)\big) \Big)=\hat{\phi}(2-pl^{(3)})=\cdots=\hat{\phi}\Big((q-1)-pl^{(q)}\Big)\\
	&\adj\hat{\phi}\Big((q-p)-1-p\big(l^{(q)}-1\big)\Big)\exch\hat{\phi}\Big(1+p(l^{(q)}-1)\Big) \quad \forall \quad l\in \mathbb{Z}.
\end{align}
Therefore by the similar argument 
as in \eqref{eq11} and by applying \textbf{Sign Rule}, equations \eqref{eq29} and \eqref{eq30} yields that
 \begin{align}\label{eq31}
 	\quad \forall \quad l\in \mathbb{Z},\quad 
 	\begin{cases}
	&\hat{\phi}(\pm 1+pl)=\hat{\phi}(\pm 2+pl)=\cdots=\hat{\phi}\big(\pm (q-1)+pl\big)=0,\\
	&\hat{\phi}(\pm q+pl)=\hat{\phi}(\pm 2q+pl)=\cdots=\hat{\phi}(\pm mq+pl)=0,\\
	&\hat{\phi}(\pm (q+1)+pl)=\hat{\phi}(\pm (2q+1)+pl)=\cdots=\hat{\phi}\Big(\pm \big((m-1)q+1\big)+pl\Big)=0.
	\end{cases}
\end{align}
Now for the remaining terms, by using \eqref{eq18} and \eqref{eq31} we have 
\begin{align}\label{eq32}
	\quad \forall \quad l\in \mathbb{Z},\quad 
	\begin{cases}
	&\hat{\phi}(q+2+pl)\adj\hat{\phi}\Big((q-p+2)+p(l+1)\Big)\exch\hat{\phi}(-2-p(l+1))=0,\\
	&\hspace*{2in}\vdots\\
	&\hat{\phi}(q+(q-1)+pl)\adj\hat{\phi}\Big(\big(q-p+(q-1)\big)+p(l+1)\Big)\\
	&\hspace{1.25in} \exch\hat{\phi}\big(-(q-1)-p(l+1)\big)=0.
	\end{cases}
\end{align}
Similarly, by employing the similar argument as in \eqref{eq32} we conclude
\begin{align}\label{eq33}
	\quad \forall \quad l\in \mathbb{Z},\quad 
	\begin{cases}
	&\hat{\phi}(2q+2+pl)=\hat{\phi}(2q+3+pl)=\cdots=\hat{\phi}(3q-1+pl)=0,\\
	&\hspace*{2in}\vdots\\
	&\hat{\phi}\big((m-2)q+2+pl\big)=\hat{\phi}\big((m-2)q+3+pl\big)=\cdots=\hat{\phi}\big((m-1)q-1+pl\big).
	\end{cases}
\end{align}
Moreover, again by applying \eqref{eq18} and \eqref{eq31} we conclude
\begin{align}\label{eq34}
	\quad \forall \quad l\in \mathbb{Z},\quad 
	\begin{cases}
	&\hat{\phi}\big((m-1)q+1+1+pl\big)=\hat{\phi}\big((p-q)+1+pl\big)\exch\hat{\phi}\big(-1-pl\big)=0,\\
	&\hspace*{2in}\vdots\\
	&\hat{\phi}\big((m-1)q+1+(q-1)+pl\big)=\hat{\phi}\big((p-q)+(q-1)+pl\big)\\
	&\hspace{2in} \exch\hat{\phi}\big(-(q-1)-pl\big)=0.
	\end{cases}
\end{align}
Finally, combining all the conditions obtained in \eqref{eq31}, \eqref{eq32}, \eqref{eq33} and \eqref{eq34} we get
\begin{align*}
	\hat{\phi}(r+pl)=0\quad \forall \quad l\in \mathbb{Z} \quad  \text{and} \quad 1\leq r\leq p-1.
\end{align*}
Summing up we have the following result.
\begin{thm}\label{mainthm}
	Let $\phi(z)=\sum\limits_{k=0}^{\infty}\hat{\phi}(z)z^k\in L^\infty(\T) $, and let  $T_\phi$ be the Toeplitz operator corresponding to the symbol $\phi$. Let $p,i,j\in \mathbb{N}$ be such that  $0\leq i<j<p$, and let $C_{p}^{i,j}$ be the corresponding conjugation defined as in \eqref{splconju} on $H^2_{\mathbb{C}}(\mathbb{D})$. If either
	\begin{enumerate}[(i)]
		\item $p$ is even and $\vert i-j\vert =\frac{p}{2}$, or
		\item $p=mq+1$ for some natural number $m\geq 2$ such that $i=q-1$ and $j=p-1$,
	\end{enumerate} 
then  $T_\phi$ is $C_{p}^{i,j}$-symmetric if and only if
\begin{align*}
\hat{\phi}(pl)&=\hat{\phi}(-pl),\quad \text{and}\quad
	\hat{\phi}(r+pl)=0\quad \forall \quad l\in \mathbb{Z} \quad  \text{and} \quad 1\leq r\leq p-1.
\end{align*}
\end{thm}
\begin{rmrk}
We expect that the above Theorem \ref{mainthm} is also valid for any $0\leq i<j<p$ and we leave this case as a work for future investigation.
\end{rmrk}

\section{Conjugation related to Model spaces}
In this section we consider a special type of conjugation on $H^2_{\mathbb{C}}(\mathbb{D})$ different from those discussed in the previous section which essentially arose from the study of natural conjugation in model spaces. Let $p=n\in \mathbb{N}$ and consider the special permutation  
\begin{center}
$\sigma:$ $\begin{pmatrix}
	a_{nk} & a_{nk+1} & \cdots & a_{nk+m} & \cdots & a_{nk+(n-1)} \\
	a_{nk+(n-1)} & a_{nk+(n-2)} & \cdots & a_{nk+(n-m-1)} & \cdots & a_{nk} 
\end{pmatrix}$	
\end{center}
on the set $\{a_{nk},a_{nk+1},\ldots,a_{nk+(n-1)}\}$ for $k\in \mathbb{N} \cup \{0\}$, and $0\leq m\leq n-1$. 
Then from \eqref{mainconju} it follows that 
\begin{equation}\label{eq35}
	C_n\equiv C_\sigma \left(\sum_{k=0}^{\infty}\sum_{m=0}^{n-1}a_{nk+m}z^{nk+m}\right)= \sum_{k=0}^{\infty}\sum_{m=0}^{n-1}\bar{a}_{nk+(n-m-1)}z^{nk+m},
\end{equation}
where $\sum\limits_{k=0}^{\infty}\sum\limits_{m=0}^{n-1}a_{nk+m}z^{nk+m}\in H^2_{\mathbb{C}}(\mathbb{D})$. As earlier,  it is easy to verify that $C_n$ is a conjugation on $H^2_{\mathbb{C}}(\mathbb{D})$. Our main aim in this section is to provide a necessary and sufficient conditions on the symbol $\phi\in L^{\infty}(\mathbb{T})$ whenever the Toeplitz operator $T_{\phi}$ is complex symmetric with respect to the conjugation $C_n$.  Let $\phi \in L^\infty(\T)$ and let $\phi (z)=\sum\limits_{k=-\infty}^{\infty}\hat{\phi}(k)z^k$. Now we assume that the Toeplitz operator $T_\phi$ is complex symmetric with respect to this conjugation $C_n$. Therefore	
\begin{align}\label{eq36}
	C_n T_\phi C_n =T^*_\phi, \quad \text{that is} \quad T_\phi C_n = C_n T_{\bar{\phi}}.
\end{align}	
It is well known that $\{z^{nk+a}:0\leq a\leq {n-1}, k\geq 0  \}$ is an orthonormal basis of the Hardy space $\hdcc$. Now applying the definition of $C_n$ (see \eqref{eq35}) it follows that	
\begin{align}\label{eq37}
	C_nT_{\bar{\phi}}(z^{nk+a}) & =C _n\left(\sum_{j=0}^{\infty}\overline{\hat{\phi}(nk+a-j)}z^j\right)
	=\sum_{j=0}^{\infty}\sum_{m=0}^{n-1}\hat{\phi}\Big(nk-nj-\{(n-1)-a\}+m\Big)z^{nj+m},
\end{align}
and	
\begin{align}\label{eq38}
	T_\phi C_n(z^{nk+a}) & = T_\phi \left(z^{nk+(n-1)-a}\right)
	=\sum_{j=0}^{\infty} \sum_{m=0}^{n-1}\hat{\phi}\Big(nj-nk-\{(n-1)-a\}+m\Big)z^{nj+m}.
\end{align}
Therefore by substituting $k-j=l\in \mathbb{Z}$ in \eqref{eq37} and \eqref{eq38}, and using \eqref{eq36} we obtain
\begin{align}\label{eq39}
	\hat{\phi}\Big(nl-\{(n-1)-a\}+m\Big)=\hat{\phi}\Big(-nl-\{(n-1)-a\}+m\Big),
\end{align}
where $0\leq a \leq {n-1}$ and $0\leq m \leq {n-1}$. In particular if $a=n-1$ and $m=0$, then the above equation \eqref{eq39} yields
\begin{equation}\label{eq40}
\hat{\phi}\Big(nl\Big) = \hat{\phi}\Big(-nl\Big) \quad \forall \quad l\in \mathbb{Z}.	
\end{equation}
The following table is essential in the sequel consisting of two different values of $a$, namely $a=0$ and $a=n-1$.
\begin{equation}\label{table}
	\begin{tabular}{ c|c } 
		\hline
		$a=0$ & $a=n-1$ \\ 
		\hline
		$\hat{\phi}\Big(nl-(n-1)\Big)$ = $\hat{\phi}\Big(-nl-(n-1)\Big)$ & $\hat{\phi}\Big(nl\Big)$ = $\hat{\phi}\Big(-nl\Big)$ \\ 
		$\hat{\phi}\Big(nl-(n-2)\Big)$ = $\hat{\phi}\Big(-nl-(n-2)\Big)$ & $\hat{\phi}\Big(nl+1\Big)$ = $\hat{\phi}\Big(-nl+1\Big)$ \\ 
		$\vdots $ &  $\vdots $\\
		$\hat{\phi}\Big(nl-(n-r)\Big)$ = $\hat{\phi}\Big(-nl-(n-r)\Big)$ & $\hat{\phi}\Big(nl+r\Big)$ = $\hat{\phi}\Big(-nl+r\Big)$ \\
		$\vdots$  & $\vdots$ \\
		$\hat{\phi}\Big(nl-1\Big)$ = $\hat{\phi}\Big(-nl-1\Big)$ & $\hat{\phi}\Big(nl+(n-2)\Big)$ = $\hat{\phi}\Big(-nl+(n-2)\Big)$ \\
		$\hat{\phi}\Big(nl\Big)$ = $\hat{\phi}\Big(-nl\Big)$ &  $\hat{\phi}\Big(nl+(n-1)\Big)$ = $\hat{\phi}\Big(-nl+(n-1)\Big)$,\\
		\hline
	\end{tabular}
\end{equation}
where $0\leq r \leq n-1$ and $l\in \mathbb{Z}$. Therefore by using the above table \eqref{table} we conclude
\begin{align}\label{eq41}
	\hat{\phi}\big(r+nl\big) = \hat{\phi}\big(r-nl\big)
	= \hat{\phi}\big(r-n -n(l-1) \big)
	= \hat{\phi}\big(r-n +n(l-1) \big)
	= \hat{\phi}\big(r+n(l-2)\big) \quad \forall \quad l\in \mathbb{Z},
\end{align}
where $1\leq r\leq n-1$. On the other hand note that $\phi \in L^\infty(\T) \subset L^2(\T)$ and hence $\sum\limits_{k=-\infty}^{\infty}\vert{\hat{\phi}(k)}\vert ^2 < \infty$. As a result the equation \eqref{eq41} yields
\begin{equation*}
\hat{\phi}\big(r+nl\big) =0 \quad \forall \quad l\in \mathbb{Z}\quad \text{and}\quad 1\leq r\leq n-1.	
\end{equation*}
So, if we assume that  $T_\phi$ is $C_n$ symmetric, then as a necessary conditions of this fact we obtain 
\begin{equation*}
	\hat{\phi}\Big(nl\Big) = \hat{\phi}\Big(-nl\Big)\quad \text{and} \quad \hat{\phi}\Big(r+nl\Big)=0, \quad \text{for any}\quad l\in \mathbb{Z}\quad \text{and}\quad 1\leq r\leq n-1.
\end{equation*}
Now conversely, if we assume $\phi \in L^\infty(\mathbb{T})$ be such that  $\hat{\phi}\Big(nl\Big) = \hat{\phi}\Big(-nl\Big)$ and $\hat{\phi}\Big(r+nl\Big)=0$ for any $ l \in \Z $ and $1\leq r \leq n-1$, Then using the definition of the conjugation $C_n$ (see \eqref{eq35}) one can easily check that
\begin{align*}
	&\Big(C_nT_{\bar{\phi}} - T_\phi C_n\Big)(z^{nk+a}) \\
	=&\sum_{j=0}^{\infty}\sum_{m=0}^{n-1}\Bigg(\hat{\phi}\Big(nk-nj-\{(n-1)-a\}+m\Big) - \hat{\phi}\Big(nj-nk-\{(n-1)-a\}+m\Big)\Bigg)z^{nj+m}=0,
\end{align*}
for any $0\leq a\leq {n-1}$ and $k\geq 0 $. Combining all the above observations we have the following main result in this section.
\begin{thm}\label{thm2}
	Let $\phi (z)=\sum\limits_{k=-\infty}^{\infty}\hat{\phi}(k)z^k \in L^\infty (\T)$, and let $T_\phi$ be the Toeplitz operator on $\hdcc$ corresponding to the symbol $\phi$. Then $T_\phi$ is complex symmetric with respect to the conjugation $C_n$ if and only if $\hat{\phi}\Big(nl\Big) = \hat{\phi}\Big(-nl\Big)~and~\hat{\phi}\Big(r+nl\Big)=0$, for any $l\in \Z$ and $1\leq r \leq n-1$.
\end{thm}

\section{Conjugations in $\hdct$}
In this section we study complex symmetric block Toeplitz operators on $\hdct$ with respect to some special conjugations on $\hdct$ introduced earlier in \cite{CamaGP,KangKoLee,LeeKo2019}. Let $C_2$ be a conjugation on $H^2_{\mathbb{C}}(\mathbb{D})$ defined as in \eqref{eq35} corresponding to $n=2$. Next we define a map  $C:\hdct \longrightarrow \hdct$ whose block matrix representation is the following:
\begin{equation}\label{conjuvec1}
C= \frac{1}{\sqrt{2}}\begin{bmatrix}
C_2 & C_2 \\
C_2 & -C_2
\end{bmatrix}
\end{equation}
Then it is important to observe that $C$ is a conjugation on $\hdct$ (see Corollary 2.8. \cite{LeeKo2019}). For more on $2\times 2$ conjugation matrices on $\mathcal{H}\oplus \mathcal{H}$ we refer to \cite{LeeKo2019}, where $\mathcal{H}$ is any complex Hilbert space. Let $\Phi \in L^\infty _{M_2}(\T)$ be such that 
$	\Phi =\begin{bmatrix}
		\phi _1 & \phi _2 \\
		\phi_3 & \phi _4
	\end{bmatrix} 
$, where $\phi_i\in L^{\infty}(\mathbb{T})$ for $1\leq i\leq 4$, and let 
$ T_\Phi= \begin{bmatrix}
	T_{\phi _1} & T_{\phi _2} \\
	T_{\phi _3} & T_{\phi _4}
\end{bmatrix}$
be corresponding block Toeplitz operator on $H^2_{\mathbb{C}^2}(\mathbb{D})$.
First we assume that the Toeplitz operator $T_\Phi$ is complex symmetric with respect to the conjugation $C$, that is
\begin{align*}
T_\Phi C & = C T_{\Phi ^*}\\
\implies \begin{bmatrix}
T_{\phi _1} & T_{\phi _2} \\
T_{\phi _3} & T_{\phi _4}
\end{bmatrix}
\frac{1}{\sqrt{2}}\begin{bmatrix}
C_2 & C_2 \\
C_2 & -C_2
\end{bmatrix} & =
\frac{1}{\sqrt{2}}\begin{bmatrix}
C_2 & C_2 \\
C_2 & -C_2
\end{bmatrix} 
\begin{bmatrix}
T_{\phi _1} & T_{\phi _2} \\
T_{\phi _3} & T_{\phi _4}
\end{bmatrix}
\end{align*}
which implies
\begin{align}\label{eq42}
C_2 (T_{\phi_1}+T_{\phi_2}) C_2  = T_{\bar{\phi}_1}+T_{\bar{\phi}_2}, \quad 
 (T_{\phi_1}-T_{\phi_2}) C_2  =C_2(T_{\bar{\phi}_3}+T_{\bar{\phi}_4}),
\end{align}
and
\begin{align}\label{eq43}
 (T_{\phi_3}+T_{\phi_4}) C_2  = C_2(T_{\bar{\phi}_1}-T_{\bar{\phi}_2}), \quad
C_2 (T_{\phi_3}-T_{\phi_4}) C_2  = T_{\bar{\phi}_3}-T_{\bar{\phi}_4}.
\end{align}
Thus by applying Theorem~\ref{thm2} for $n=2$ and using equations \eqref{eq42} and \eqref{eq43} we conclude
\vspace{0.1in}

$(i)$ $T_{\phi _1 +\phi _2}$ is $C_2$ symmetric, that is $\widehat{\phi_1+\phi_2}(2l)=\widehat{\phi_1+\phi_2}(-2l)$ and $\widehat{\phi_1+\phi_2}(2l+1)=0$, for all $l\in \Z$,
\vspace{0.1in}

$(ii)$ $T_{\phi _3 -\phi _4}$ is $C_2$ symmetric, that is $\widehat{\phi_3-\phi_4}(2l)=\widehat{\phi_3-\phi_4}(-2l)$ and $\widehat{\phi_3-\phi_4}(2l+1)=0$, for all $l\in \Z$,
\vspace{0.1in}

$(iii)$ $C_2(T_{\phi_1}-T_{\phi_2}) C_2 =T_{\bar{\phi}_3}+T_{\bar{\phi}_4}$.
\vspace{0.1in}

\noindent Suppose $\phi _1 -\phi_2 =\psi _1$ and $\phi _3 +\phi_4 =\psi _2$, and let $\psi_1(z)=\sum\limits _{n=-\infty}^{\infty}\hat{\psi}_1(n)z^n$ and $ \psi_2(z)=\sum\limits _{n=-\infty}^{\infty}\hat{\psi}_2(n)z^n $.
Then the above condition $(iii)$ becomes
\begin{align}\label{eq44}
T_{\psi _1}C_2 & = C_2 T_{\bar{\psi}_2}.
\end{align}
Next by applying the definition of $C_2$ and using \eqref{eq44} we get for any $m\geq 0$ that
\begin{align*}
T_{\psi _1}C_2(z^{2m}) & = C_2 T_{\bar{\psi}_2}(z^{2m})\\
\Rightarrow T_{\psi _1}\left(z^{2m+1}\right)& = C_2\left(\sum_{n=0}^{\infty}\overline{\hat{\psi}_2(-(n-2m))}z^n \right)\\
\Rightarrow \sum_{n=0}^{\infty}\hat{\psi}_1(n-2m-1)z^n & = \sum_{k=0}^{\infty}\hat{\psi}_2(-(2k-2m))z^{2k+1} +\sum_{k=0}^{\infty}\hat{\psi}_2(-(2k+1-2m))z^{2k},
\underline{\textit{}}\end{align*}
which by equating the Fourier coefficient yields the following conditions:
\begin{align}\label{eq45}
	\forall \quad k,m\geq 0,\quad 
\begin{cases}
	& \hat{\psi}_1(2k-2m-1) = \hat{\psi}_2(-(2k+1-2m)),\\
	& \hat{\psi}_1(2k-2m) = \hat{\psi}_2(-(2k-2m)).
\end{cases}
\end{align}
Similarly, for any $m\geq 0$ we also get
\begin{align*}
T_{\psi _1}C_2(z^{2m+1}) & = C_2 T_{\bar{\psi}_2}(z^{2m+1})\\
\Rightarrow T_{\psi _1}\left(z^{2m}\right)& = C_2\left(\sum_{n=0}^{\infty}\overline{\hat{\psi}_2(-(n-2m-1))}z^n \right)\\
\Rightarrow \sum_{n=0}^{\infty}\hat{\psi}_1(n-2m)z^n & = \sum_{k=0}^{\infty}\hat{\psi}_2(-(2k-2m-1))z^{2k+1} +\sum_{k=0}^{\infty}\hat{\psi}_2(-(2k+1-2m-1))z^{2k},
\end{align*}
which leads to the following conditions:
\begin{align}\label{eq46}
	\forall \quad k,m\geq 0,\quad 
	\begin{cases}
		& \hat{\psi}_1(2k-2m+1) = \hat{\psi}_2(-(2k-2m-1)),\\
		& \hat{\psi}_1(2k-2m) = \hat{\psi}_2(-(2k-2m)).
	\end{cases}
\end{align}
Substituting $k-m=l \in \Z$ in \eqref{eq45} and \eqref{eq46} we have the following set of conditions 
\begin{align}\label{eq47}
	\forall \quad l\in \mathbb{Z},\quad 
	\begin{cases}
		& \widehat{\phi_1-\phi_2}(2l) =\widehat{\phi_3+\phi_4}(-2l),\\
		& \widehat{\phi_1-\phi_2}(2l-1) =\widehat{\phi_3+\phi_4}(-2l-1),\\
		& \widehat{\phi_1-\phi_2}(2l+1) =\widehat{\phi_3+\phi_4}(-2l+1).
	\end{cases}
\end{align}
Moreover, using the above mentioned conditions in \eqref{eq47} we conclude 
\begin{align*}
\widehat{\phi_1-\phi_2}(2l+1)& =\widehat{\phi_3+\phi_4}(-2l+1)
=\widehat{\phi_3+\phi_4}\left(-2(l-1)-1\right)\\
& =\widehat{\phi_1-\phi_2}\left(2(l-1)-1\right)
=\widehat{\phi_1-\phi_2}\left(2(l-2)+1\right)\quad \forall \quad l\in \mathbb{Z},
\end{align*}
and hence $\widehat{\phi_1-\phi_2}(2l+1)=0$, for all $l\in \mathbb{Z}$ since $\phi_1-\phi_2\in L^{\infty}(\mathbb{T})\subseteq L^2(\mathbb{T})$. Similarly, we also conclude $\widehat{\phi_3+\phi_4}(2l+1)=0$, for all $l\in \Z$.
Consequently, combining all the above obtained conditions we get
\begin{align*}
	\forall \quad l\in \mathbb{Z},\quad 
	\begin{cases}
		& \widehat{\phi_1+\phi_2}(2l)=\widehat{\phi_1+\phi_2}(-2l), \quad \widehat{\phi_1+\phi_2}(2l+1)=0,\\
		& \widehat{\phi_3-\phi_4}(2l)=\widehat{\phi_3-\phi_4}(-2l), \quad \widehat{\phi_3-\phi_4}(2l+1)=0,\\
		& \widehat{\phi_1-\phi_2}(2l)=\widehat{\phi_3+\phi_4}(-2l),\quad  \widehat{\phi_1-\phi_2}(2l+1)=0,\quad  \widehat{\phi_3+\phi_4}(2l+1)=0,
	\end{cases}
\end{align*}
which after slight modifications reduces to,
\begin{align}\label{eq48}
	\forall \quad l\in \mathbb{Z},\quad 
	\begin{cases}
		& \widehat{\phi_1+\phi_2}(2l)=\widehat{\phi_1+\phi_2}(-2l),\quad \widehat{\phi_3-\phi_4}(2l)=\widehat{\phi_3-\phi_4}(-2l),\\
		& \widehat{\phi_1-\phi_2}(2l)=\widehat{\phi_3+\phi_4}(-2l),\quad \hat{\phi}_i(2l+1)=0\quad \text{for} \quad 1\leq i\leq 4.
	\end{cases}
\end{align}
Summing up, we have the following theorem.
\begin{thm}\label{thm3}
	Let 
	$\Phi =\begin{bmatrix}
	\phi _1 & \phi _2 \\
	\phi_3 & \phi _4
	\end{bmatrix} \in L^\infty _{M_2}(\T)$, and let $C:H^2_{\mathbb{C}^2}(\mathbb{D})\longrightarrow H^2_{\mathbb{C}^2}(\mathbb{D})$ be a conjugation on $H^2_{\mathbb{C}^2}(\mathbb{D})$ whose block matrix representation is
	$\frac{1}{\sqrt{2}}\begin{bmatrix}
	C_2 & C_2 \\
	C_2 & -C_2
	\end{bmatrix}$, where $C_2$ is conjugation on $H^2_{\mathbb{C}}(\mathbb{D})$ defined as in \eqref{eq35}. Then the Toeplitz operator $T_\Phi$ is  complex symmetric with respect to the conjugation $C$ if and only if 
\begin{align*}
&\widehat{\phi_1+\phi_2}(2l)=\widehat{\phi_1+\phi_2}(-2l),\quad \widehat{\phi_3-\phi_4}(2l)=\widehat{\phi_3-\phi_4}(-2l)\quad\widehat{\phi_1-\phi_2}(2l)=\widehat{\phi_3+\phi_4}(-2l),\\
&\hat{\phi}_i(2l+1)=0, \text{ for all } l\in \Z \text{ and } 1\leq i\leq 4.
\end{align*}
\end{thm}
Alternatively, by adding and subtracting the conditions obtained in \eqref{eq48} we have the following theorem. 
\begin{thm}\label{thm4}
	Let 
	$\Phi =\begin{bmatrix}
		\phi _1 & \phi _2 \\
		\phi_3 & \phi _4
	\end{bmatrix} \in L^\infty _{M_2}(\T)$, and let $C:H^2_{\mathbb{C}^2}(\mathbb{D})\longrightarrow H^2_{\mathbb{C}^2}(\mathbb{D})$ be a conjugation on $H^2_{\mathbb{C}^2}(\mathbb{D})$ whose block matrix representation is
	$\frac{1}{\sqrt{2}}\begin{bmatrix}
		C_2 & C_2 \\
		C_2 & -C_2
	\end{bmatrix}$,  where $C_2$ is conjugation on $H^2_{\mathbb{C}}(\mathbb{D})$ defined as in \eqref{eq35}. Then the Toeplitz operator $T_\Phi$ is  complex symmetric with respect to the conjugation $C$ if and only if 
	\begin{align*}
		&\widehat{\phi_1+\phi_2}(2l)=\widehat{\phi_1+\phi_2}(-2l),\quad \widehat{\phi_1+\phi_3}(2l)=\widehat{\phi_1+\phi_3}(-2l)\quad\widehat{\phi_1+\phi_4}(2l)=\widehat{\phi_1+\phi_4}(-2l),\\
		&\hat{\phi}_i(2l+1)=0, \text{ for all } l\in \Z \text{ and } 1\leq i\leq 4.
	\end{align*}
\end{thm}
 \begin{rmrk}
 	Let 
 $\Phi =\begin{bmatrix}
 	\phi _1 & \phi _2 \\
 	\phi_3 & \phi _4
 \end{bmatrix} \in L^\infty _{M_2}(\T)$, and let 
$ T_\Phi= \begin{bmatrix}
	T_{\phi _1} & T_{\phi _2} \\
	T_{\phi _3} & T_{\phi _4}
\end{bmatrix}$
be the corresponding block Toeplitz operator on $H^2_{\mathbb{C}^2}(\mathbb{D})$. Then using Theorem~\ref{thm4} and Theorem~\ref{thm2} we conclude the following: If $T_{\phi _i}$ is complex symmetric with respect to the conjugation $C_2$ on  $H^2_{\mathbb{C}}(\mathbb{D})$ for all $1\leq i\leq 4$, then $T_\Phi$ is complex symmetric with respect to the conjugation $C$ on  $H^2_{\mathbb{C}^2}(\mathbb{D})$. Conversely,  if  $T_\Phi$ is complex symmetric with respect to the conjugation $C$ on  $H^2_{\mathbb{C}^2}(\mathbb{D})$ and if one of $T_{\phi_i}$ is complex symmetric with respect to the conjugation $C_2$ on  $H^2_{\mathbb{C}}(\mathbb{D})$, then  rest of $T_{\phi_i}$ is also complex symmetric with respect to the conjugation $C_2$ on  $H^2_{\mathbb{C}}(\mathbb{D})$.
\end{rmrk}
Next, let $C_1$ and $C_2$ be two conjugations on $H^2_{\mathbb{C}}(\mathbb{D})$ defined as in \eqref{eq35} corresponding to $n=1$ and $n=2$ respectively. 
Then it is easy to verify that $C_1$ commutes with $C_2$, that is $C_1C_2= C_2C_1$.
Let $\widetilde{C}:\hdct \to \hdct$ be a conjugation on $\hdct$ defined in \cite[Corollary 2.8]{LeeKo2019} whose block matrix representation is the following
\begin{align}\label{conjuvec2}
\widetilde{C} & = \frac{1}{\sqrt{2}}\begin{bmatrix}
C_2 & C_1 \\
C_1 & -C_2
\end{bmatrix}.
\end{align}
Our next aim is to investigate the complex-symmetry of the Toeplitz operator $T_\Phi :\hdct \to \hdct$ having symbol 
$\Phi =\begin{bmatrix}
\phi _1 & \phi _2 \\
\phi_3 & \phi _4
\end{bmatrix} \in L^\infty _{M_2}(\T)$ with respect to the conjugation $\widetilde{C}$.
As earlier, we first assume $T_\Phi$ is complex-symmetric with respect to $\widetilde{C}$, that is
\begin{align*}
T_\Phi \widetilde{C} &= \widetilde{C} T_{\Phi^*}\\
\implies\begin{bmatrix}
T_{\phi _1} & T_{\phi _2} \\
T_{\phi _3} & T_{\phi _4}
\end{bmatrix}\begin{bmatrix}
C_2 & C_1 \\
C_1 & -C_2
\end{bmatrix} &= \begin{bmatrix}
C_2 & C_1 \\
C_1 & -C_2
\end{bmatrix}\begin{bmatrix}
T_{\bar{\phi _1}} & T_{\bar{\phi _3}} \\
T_{\bar{\phi _2}} & T_{\bar{\phi _4}}
\end{bmatrix}
\end{align*} 
which yields the following set of conditions:
\begin{equation}\label{eq49}
\begin{cases}
T_{\phi_1}C_2+T_{\phi_2}C_1 = C_2T_{\bar{\phi}_1}+C_1T_{\bar{\phi}_2},\\
T_{\phi_1}C_1-T_{\phi_2}C_2 = C_2T_{\bar{\phi}_3}+C_1T_{\bar{\phi}_4},
\end{cases}
\begin{cases}
T_{\phi_3}C_2+T_{\phi_4}C_1 = C_1T_{\bar{\phi}_1}-C_2T_{\bar{\phi}_2},\\
T_{\phi_3}C_1-T_{\phi_4}C_2 = C_1T_{\bar{\phi}_3}-C_2T_{\bar{\phi}_4}.
\end{cases}
\end{equation}
Therefore, as earlier using the definition of $C_1$, $C_2$ and using \eqref{eq49} we get for any $m\geq 0$ that
\begin{align*}
\left(T_{\phi_1}C_2+T_{\phi_2}C_1\right)(z^{2m})  = & \left(C_2T_{\bar{\phi}_1}+C_1T_{\bar{\phi}_2} \right)(z^{2m})\\
\implies \sum_{n=0}^{\infty} \hat{\phi}_1\left(n-2m-1\right)z^n + \sum_{n=0}^{\infty} \hat{\phi}_2\left(n-2m\right)z^n  = &~ C_2 \left(\sum_{n=0}^{\infty}\overline{\hat{\phi}_1\left(-(n-2m)\right)} z^n\right) \\
& \quad \quad \quad \quad + C_1 \left(\sum_{n=0}^{\infty}\overline{\hat{\phi}_2\left(-(n-2m)\right)} z^n\right)
\end{align*}
\begin{align*}
& \hspace{0.2in} \implies \sum_{k=0}^{\infty}\left\{\hat{\phi}_1\left(2k-2m\right)+ \hat{\phi}_2\left(2k+1-2m\right)\right\}z^{2k+1}  +\sum_{k=0}^{\infty}\left\{\hat{\phi}_1\left(2k-2m-1\right)+ \hat{\phi}_2\left(2k-2m\right)\right\}z^{2k} \\
& \hspace{1in} = \sum_{k=0}^{\infty}\hat{\phi}_1\left(-(2k-2m)\right)z^{2k+1} +
\sum_{k=0}^{\infty}\hat{\phi}_1\left(-(2k+1-2m)\right)z^{2k}\\
& \hspace{2in} +\sum_{k=0}^{\infty} \hat{\phi}_2\left(-(2k-2m)\right)z^{2k} +\sum_{k=0}^{\infty} \hat{\phi}_2\left(-(2k+1-2m)\right)z^{2k+1},
\end{align*}
which by substituting the index $k-m=l \in \mathbb{Z}$ and equating the Fourier coefficient leads to the following conditions:
\begin{align*}
	\begin{cases}
\hat{\phi}_1(2l) + \hat{\phi}_2(2l+1) & = \hat{\phi}_1(-2l) + \hat{\phi}_2(-2l-1),\\
\hat{\phi}_1(2l-1) + \hat{\phi}_2(2l) & = \hat{\phi}_1(-2l-1) + \hat{\phi}_2(-2l).
\end{cases}
\end{align*} 
Similarly, by repeating the above similar calculations and using the equation
 $\left(T_{\phi_1}C_2+T_{\phi_2}C_1\right)(z^{2m+1})$  $=\left(C_2T_{\bar{\phi}_1}+C_1T_{\bar{\phi}_2} \right)(z^{2m+1})$ for any $m\geq 0$ we get the following conditions:
\begin{align*}
		\forall \quad l\in \mathbb{Z},\quad
	\begin{cases}
\hat{\phi}_1(2l) + \hat{\phi}_2(2l-1) & = \hat{\phi}_1(-2l) + \hat{\phi}_2(-2l+1),\\
\hat{\phi}_1(2l+1) + \hat{\phi}_2(2l) & = \hat{\phi}_1(-2l+1) + \hat{\phi}_2(-2l).
\end{cases}
\end{align*} 
 Therefore continuing the above process, the equation $T_{\phi_1}C_1-T_{\phi_2}C_2 = C_2T_{\bar{\phi}_3}+C_1T_{\bar{\phi}_4}$ yields the following set of conditions:
 \begin{equation}\label{eq50}
 	\forall \quad l\in \mathbb{Z},\quad
 \begin{cases}
 \hat{\phi}_1(2l) - \hat{\phi}_2(2l-1)  = \hat{\phi}_4(-2l) + \hat{\phi}_3(-2l-1),\\
 \hat{\phi}_1(2l+1) - \hat{\phi}_2(2l)  = \hat{\phi}_4(-2l-1) + \hat{\phi}_3(-2l),\\
 \hat{\phi}_1(2l-1) - \hat{\phi}_2(2l)  = \hat{\phi}_4(-2l+1) + \hat{\phi}_3(-2l),\\
 \hat{\phi}_1(2l) - \hat{\phi}_2(2l+1)  = \hat{\phi}_4(-2l) + \hat{\phi}_3(-2l+1).
 \end{cases}
 \end{equation} 
Similarly, the equation $T_{\phi_3}C_2+T_{\phi_4}C_1 = C_1T_{\bar{\phi}_1}-C_2T_{\bar{\phi}_2}$
leads to the following set of conditions:
\begin{equation}\label{eq51}
	\forall \quad l\in \mathbb{Z},\quad
\begin{cases}
\hat{\phi}_4(2l) + \hat{\phi}_3(2l-1)  = \hat{\phi}_1(-2l) - \hat{\phi}_2(-2l-1),\\
\hat{\phi}_4(2l+1) + \hat{\phi}_3(2l)  = \hat{\phi}_1(-2l-1) - \hat{\phi}_2(-2l),\\
\hat{\phi}_4(2l-1) + \hat{\phi}_3(2l)  = \hat{\phi}_1(-2l+1) - \hat{\phi}_2(-2l),\\
\hat{\phi}_4(2l) + \hat{\phi}_3(2l+1)  = \hat{\phi}_1(-2l) - \hat{\phi}_2(-2l+1).
\end{cases}
\end{equation} 
Furthermore, the equation $T_{\phi_3}C_1-T_{\phi_4}C_2 = C_1T_{\bar{\phi}_3}-C_2T_{\bar{\phi}_4}$
gives the following set of conditions:
\begin{equation*}
		\forall \quad l\in \mathbb{Z},\quad
\begin{cases}
\hat{\phi}_3(2l) - \hat{\phi}_4(2l-1)  = \hat{\phi}_3(-2l) - \hat{\phi}_4(-2l-1),\\
\hat{\phi}_3(2l+1) - \hat{\phi}_4(2l)  = \hat{\phi}_3(-2l-1) - \hat{\phi}_4(-2l),\\
\hat{\phi}_3(2l-1) - \hat{\phi}_4(2l)  = \hat{\phi}_3(-2l+1) - \hat{\phi}_4(-2l),\\
\hat{\phi}_3(2l) - \hat{\phi}_4(2l+1)  = \hat{\phi}_3(-2l) - \hat{\phi}_4(-2l+1).
\end{cases}
\end{equation*} 
It is important to observe that some repetition is there in the set of conditions obtained in \eqref{eq50} and \eqref{eq51}. Thus by removing those repetition we get the following complete list of conditions:
\begin{equation*}
	\forall \quad l\in \mathbb{Z},\quad
\begin{cases}
\hat{\phi}_1(2l) + \hat{\phi}_2(2l+1)  = \hat{\phi}_1(-2l) + \hat{\phi}_2(-2l-1),\\
\hat{\phi}_1(2l-1) + \hat{\phi}_2(2l)  = \hat{\phi}_1(-2l-1) + \hat{\phi}_2(-2l),\\
\hat{\phi}_1(2l) + \hat{\phi}_2(2l-1)  = \hat{\phi}_1(-2l) + \hat{\phi}_2(-2l+1),\\
\hat{\phi}_1(2l+1) + \hat{\phi}_2(2l)  = \hat{\phi}_1(-2l+1) + \hat{\phi}_2(-2l),\\
\hat{\phi}_1(2l) - \hat{\phi}_2(2l-1)  = \hat{\phi}_4(-2l) + \hat{\phi}_3(-2l-1),\\
\hat{\phi}_1(2l+1) - \hat{\phi}_2(2l)  = \hat{\phi}_4(-2l-1) + \hat{\phi}_3(-2l),\\
\hat{\phi}_1(2l-1) - \hat{\phi}_2(2l)  = \hat{\phi}_4(-2l+1) + \hat{\phi}_3(-2l),\\
\hat{\phi}_1(2l) - \hat{\phi}_2(2l+1)  = \hat{\phi}_4(-2l) + \hat{\phi}_3(-2l+1),\\
\hat{\phi}_3(2l) - \hat{\phi}_4(2l-1)  = \hat{\phi}_3(-2l) - \hat{\phi}_4(-2l-1),\\
\hat{\phi}_3(2l+1) - \hat{\phi}_4(2l)  = \hat{\phi}_3(-2l-1) - \hat{\phi}_4(-2l),\\
\hat{\phi}_3(2l-1) - \hat{\phi}_4(2l)  = \hat{\phi}_3(-2l+1) - \hat{\phi}_4(-2l),\\
\hat{\phi}_3(2l) - \hat{\phi}_4(2l+1)  = \hat{\phi}_3(-2l) - \hat{\phi}_4(-2l+1),
\end{cases}
\end{equation*} 
which after simplifying again we obtain the following set of minimal conditions.
\begin{align*}
		\forall \quad l\in \mathbb{Z},\quad
	\begin{cases}
\widehat{\phi_1 +\phi _4}(2l+1)  = 0, \widehat{\phi_2 -\phi _3}(2l+1)  = 0, \\
\hat{\phi}_1(2l) -\hat{\phi}_1(2l+2)  = \hat{\phi}_1(-2l) -\hat{\phi}_1(-2l-2),\\
\hat{\phi}_2(2l) -\hat{\phi}_1(2l-2)  = \hat{\phi}_2(-2l) -\hat{\phi}_2(-2l+2),\\
\hat{\phi}_3(2l) -\hat{\phi}_1(2l-2)  = \hat{\phi}_3(-2l) -\hat{\phi}_3(-2l+2),\\
\hat{\phi}_4(2l) -\hat{\phi}_1(2l+2)  = \hat{\phi}_4(-2l) -\hat{\phi}_4(-2l-2),\\
\widehat{\phi_2 +\phi _3}(2l)+\widehat{\phi_1 -\phi _4}(2l-1)  = \widehat{\phi_2 +\phi _3}(-2l)+\widehat{\phi_1 -\phi _4}(-2l-1),\\
\widehat{\phi_1 -\phi _4}(2l)+\widehat{\phi_2 +\phi _3}(2l+1)  = \widehat{\phi_1 -\phi _4}(-2l)+ \widehat{\phi_2 +\phi _3}(-2l-1).
\end{cases}
 \end{align*}
Summing up, we have the following theorem in this section.
\begin{thm}\label{thm5}
	Let 
	$\Phi =\begin{bmatrix}
		\phi _1 & \phi _2 \\
		\phi_3 & \phi _4
	\end{bmatrix} \in L^\infty _{M_2}(\T)$, and let $\widetilde{C}:H^2_{\mathbb{C}^2}(\mathbb{D})\longrightarrow H^2_{\mathbb{C}^2}(\mathbb{D})$ be a conjugation on $H^2_{\mathbb{C}^2}(\mathbb{D})$ whose block matrix representation is
	$\frac{1}{\sqrt{2}}\begin{bmatrix}
		C_2 & C_1 \\
		C_1 & -C_2
	\end{bmatrix}$, where $C_1$ and $C_2$ are conjugations on $H^2_{\mathbb{C}}(\mathbb{D})$ defined as in \eqref{eq35}. Now if the Toeplitz operator $T_\Phi$ is  complex symmetric with respect to the conjugation $\widetilde{C}$, then for any $l\in \mathbb{Z}$ we get
\begin{align*}
	\widehat{\phi_1 +\phi _4}(2l+1)  & = 0, \widehat{\phi_2 -\phi _3}(2l+1)  = 0, \\
	\hat{\phi}_1(2l) -\hat{\phi}_1(2l+2) & = \hat{\phi}_1(-2l) -\hat{\phi}_1(-2l-2),\\
	\hat{\phi}_2(2l) -\hat{\phi}_1(2l-2) & = \hat{\phi}_2(-2l) -\hat{\phi}_2(-2l+2),\\
	\hat{\phi}_3(2l) -\hat{\phi}_1(2l-2) & = \hat{\phi}_3(-2l) -\hat{\phi}_3(-2l+2),\\
	\hat{\phi}_4(2l) -\hat{\phi}_1(2l+2) & = \hat{\phi}_4(-2l) -\hat{\phi}_4(-2l-2),\\
	\widehat{\phi_2 +\phi _3}(2l)+\widehat{\phi_1 -\phi _4}(2l-1) & = \widehat{\phi_2 +\phi _3}(-2l)+\widehat{\phi_1 -\phi _4}(-2l-1),\\
	\widehat{\phi_1 -\phi _4}(2l)+\widehat{\phi_2 +\phi _3}(2l+1) & = \widehat{\phi_1 -\phi _4}(-2l)+ \widehat{\phi_2 +\phi _3}(-2l-1).
	\end{align*}
\end{thm}

\section{Concluding Remarks}
It is important to observe that in section 2 and section 3,  we give a characterization of complex symmetric Toeplitz operator $T_{\phi}$ on $H^2_{\mathbb{C}}(\mathbb{D})$ with respect to a class of conjugations that we obtain as a special cases of $C_{\sigma}$ defined in \eqref{mainconju}. Therefore it is natural to ask the following question in the sequel:
\begin{align*}
\textbf{\text{Question:}}~ &\text{Characterize a complex symmetric Toeplitz operator $T_{\phi}$ on the Hardy space}~ H^2_{\mathbb{C}}(\mathbb{D})~\\
&\text{with respect to the conjugation}~C_{\sigma}~\text{defined in}~ \eqref{mainconju}
\end{align*}

\noindent We expect to have similar type of characterizations as obtained in Theorem~\ref{mainthm}, Theorem~\ref{thm2} and leave this as a subject for future investigation.

	\section*{Acknowledgements}
	\textit{The research of the first named author is supported by the Mathematical Research Impact Centric Support (MATRICS) grant, File No :MTR/2019/000640, by the Science and Engineering Research Board (SERB), Department of Science $\&$ Technology (DST), Government of India. The second and the third named author gratefully acknowledge the support provided by IIT Guwahati, Government of India. The research of the fourth named author is
supported by DST-INSPIRE Faculty Fellowship No. - DST/INSPIRE/04/2019/000769.}


\begin{thebibliography}{99}
		
	
	\bibitem{BE07}
	 E.L.~Basor and T.~Ehrhardt, {\emph{Torsten Asymptotic of block Toeplitz determinants and the classical Dimer model}}, Comm Math Phys., 274 (2007), 427--455.
	
	
	\bibitem{BFGJ}
	C.~Bender, A.~Fring, U.~G$\ddot{u}$nther and H.~Jones, {\emph{Quantum physics with non-Hermitian operators}}, J.Phys.A, 45 (2012), 440301.
	
	\bibitem{BH}
	A.~Brown and P.R.~Halmos, {\emph{Algebraic properties of Toeplitz operators}}, J.Reine Angew. Math., 213 (1963–1964), 89--102.
	
	\bibitem{CamaGP}
	M.C.~C$\hat{a}$mara, K.~ Kli$\acute{s}$-Garlicka, and M.~Ptak, {\emph{Complex symmetric completions of partial operator matrices}}, Linear and Multilinear Algebra, 2019, DOI: 10.1080/03081087.2019.1631246.
	
	\bibitem{CurtHL}
	R. E.~Curto, I. S.~Hwang,and  W. Y.~Lee, {\emph{Which subnormal Toeplitz operators are either normal or analytic}}, J. Func. Anal., 263 (2012), 2333--2354.
	
	\bibitem{FLee}
	D.R.~Farenick and W.Y.~Lee, {\emph{Hyponormality and spectra of Toeplitz operators}}, Trans. Amer. Math. Soc., 384 (1996), 4153--4174.
	
	
	\bibitem{G06}
	S.R.~Garcia, {\emph{Conjugation and Clark operators}}, in: Contemp. Math., vol.393, 2006, pp.67--112.
	
	\bibitem{GPP14}
	S.R.~Garcia, E.~Prodan and M. ~Putinar, {\emph{Mathematical and physical aspects of complex symmetric operators}}, J.Phys.A, 47 (2014), 1--51.
	
	\bibitem{GP06}
	S.R.~Garcia and M.~Putinar, {\emph{Complex symmetric operators and applications}}, Trans. Amer. Math. Soc., 358 (2006), 1285--1315.
	
	\bibitem{GP07}
	S.R.~Garcia and M.~Putinar, {\emph{Complex symmetric operators and applications II}}, Trans. Amer. Math. Soc., 359 (2007), 3913--3931.
	
	\bibitem{GW09}
	S. R.~Garcia and W. R.~Wogen, {\emph{Complex symmetric partial isometries}}, J. Funct. Anal., 257 (2009), 1251--1260.
	
	\bibitem{GW10}
	S. R.~Garcia and W. R.~Wogen, {\emph{Some new classes of complex symmetric operators}}, Trans. Amer. Math. Soc., 362 (2010), 6065--6077.
	
	
	\bibitem{GZ}
	K.~Guo and S.~Zhu, {\emph{A canonical decomposition of complex symmetric operators}}, J.Operator Theory, 72 (2014), 529--547.
	
	
	\bibitem{KangKoLee}	
	D.~Kang, E.~Ko and J.E.~Lee, {\emph {Remarks on complex symmetric Toeplitz operators}}, Linear Multilinear Algebra, 2020, DOI: 10.1080/03081087.2020.1842847.
	
	\bibitem{LeeKo2016}	
	E.~Ko and J.E.~Lee, {\emph{On complex symmetric Toeplitz operators}}, J Math Anal Appl, 2016 (434), 20--34.
	
	\bibitem{LeeKo2019}	
	E.~Ko and J.E.~Lee, {\emph {Remark on complex symmetric operator matrices}}, Linear Multilinear Algebra, 2019 67(6),1198--1216.
	
	\bibitem{LLC}
	A.~Li, Y.~Liu, and Y.~ Chen, {\emph{Complex symmetric Toeplitz operators on the Dirichlet space}}, J Math Anal Appl, 2020 (487), 123998. 
	
	\bibitem{LYL}
	R.~ Li, Y.~Yang and Y.~ Lu, {\em {A class of complex symmetric Toeplitz operators on Hardy and Bergman spaces}}, J Math Anal Appl, 2020 (489), 124173. 
	
	\bibitem{Prunel}	
	E.~de Prunel$\acute{e}$, {\emph {Conditions for bound states in a periodic linear chain, and the spectra of a class of Toeplitz operators interms of polylogarithm functions}}, J. Phys. A 36 (2003) 8797--8815.
	
	\bibitem{Noor}
	S.~ Waleed Noor, {\emph{Complex symmetry of Toeplitz operators with continuous symbols}}, Arch. Math. 109 (2017) 455--460.
	
	
\end{thebibliography}
\end{document}